\documentclass[12pt]{amsart}
\usepackage{amssymb}
\newtheorem{thm}{Theorem}[section]

\newtheorem{lem}[thm]{Lemma}

\numberwithin{equation}{section}
\begin{document}
\title[]{On the automorphism group of a possible symmetric $(81,16,3)$ design}%
\author[A. Abdollahi et al.]{A. Abdollahi$^*$, H. R. Maimani and R. Torabi}%
\address{Department of Mathematics, University of Isfahan, Isfahan 81746-71441, Iran;
and Institute for Studies in Theoretical Physics and Mathematics
(IPM); Tehran, Iran.}%
\address{Deaprtment of Mathematics, University of Tehran,  Tehran, Iran.}
\address{Deaprtment of Mathematics, University of Tehran,  Tehran, Iran.}
\email{a.abdollahi@math.ui.ac.ir} \email{maimani@ipm.ir}%
\thanks{$^*$ Corresponding Author.}
\thanks{The research of A. Abdollahi was in part supported by a grant from IPM (No. 85200032). He also
thanks the Center of Excellence for Mathematics, University of
Isfahan.}
%\date{}%
%\dedicatory{}%
%\commby{}%
% ----------------------------------------------------------------
\begin{abstract}
In this paper we study the automorphism group of a possible
symmetric $(81,16,3)$ design.
\end{abstract}
\maketitle
% ----------------------------------------------------------------
\section{Introduction}
Let $v$, $k$ and $\lambda$ be non-negative integers such that
$v>k>\lambda$. By a symmetric $(v,k,\lambda)$ design, we mean a
pair $D=(V,\mathcal{B})$, where  $V$ is a $v$-set and
$\mathcal{B}$ is a set of $k$-subsets of $V$ such that the
following four requirements are satisfied by $D$:
\begin{enumerate}
\item $|\mathcal{B}|=v$.
\item any element of $V$ belongs to  precisely $k$ members of
$\mathcal{B}$.
\item any two distinct members of $\mathcal{B}$ intersect in
exactly $\lambda$ elements of $V$.
\item any two distinct elements  of $V$ are in
exactly $\lambda$ members  of $\mathcal{B}$.
\end{enumerate}
As usual, the elements of $V$ are called points of $D$ and the
members of $\mathcal{B}$ are called  blocks of the design $D$. An
automorphism of a symmetric  design $D=(V,\mathcal{B})$ is a
permutation on $V$ which sends blocks to blocks. The set of all
automorphisms of $D$ with the composition rule of maps forms the
full automorphism group of $D$ which will be denoted by $Aut(D)$.
 If $\alpha$ is an
automorphism of $D$, we denote by $F(\alpha)$ the set of all
points which are fixed by $\alpha$; and $F_b(\alpha)$ denotes the
set of all blocks which are fixed by $\alpha$.\\

 Over the years, researchers have tackled problems related to
symmetric designs. The question of existence still remains
unsettled for many parameter sets. Indeed, if we list the
parameters $(v,k,\lambda)$ in order of increasing $n=k-\lambda$,
then $(81,16,3)$ would be the smallest unknown case \cite{Trung}.
On the other hand, the success of almost all the design
construction methods depends  heavily on a proper choice of
possible automorphism groups \cite{Camina}.\\

As far as we know, the only known results on a possible
$(81,16,3)$ design are  the following:
\begin{thm}
{\rm (See \cite{Arasu})} There is no symmetric $(81,16,3)$ design
with an abelian regular $3$-group of automorphisms.
\end{thm}
\begin{thm}\label{fix2}{\rm (See \cite{Marangunic})}
Let $\alpha$ be an automorphism  of a possible symmetric
$(81,16,3)$ design of order $2$. Then $|F(\alpha)|=9$.
\end{thm}
\begin{thm}{\rm (See \cite{Garapic})} The alternating group $A_5$
of degree $5$ cannot be isomorphic to a group of automorphisms  of
a possible symmetric $(81,16,3)$ design.
\end{thm}
T. Spence has announced in his home page
  \begin{center} {\tt http://www.maths.gla.ac.uk/\~es/}
  \end{center}
 that there is no symmetric
$(81,16,3)$ designs having a ``certain'' fixed-point free
automorphism of order 3.\\
Our main  result is:
\begin{thm}\label{mainmain}
If $G$ is the full automorphism group of  a possible symmetric
$(81,16,3)$ design, then $|G|=2^{\alpha} 3^{\beta} 5^{\gamma}
13^{\sigma}$, where $ \gamma\leq 1$, $\sigma \leq 1$. Moreover,
$G$ has no subgroup of order  $65$, and has no elements of orders
$10$ or $26$; and $G$ does not contain any abelian $2$-subgroup of
rank greater than $3$.
\end{thm}
In Section 2, some general results on the automorphism groups of
a symmetric design are given and in Section 3, we prove a series
of Lemmata.  Based on them we can prove Theorem \ref{mainmain}.
\section{Some general results on the automorphism group of a  symmetric  design }
\begin{lem}{\rm (See \cite{Lander})} \label{fix=fix} Let $\alpha$ be an automorphism of a nontrivial symmetric
$(v,k,\lambda)$ design. Then $|F(\alpha)|=|F_b(\alpha)|$.
\end{lem}
\begin{lem}\label{fix}{\rm (See \cite[Corollary 3.7, p. 82]{Lander})}
Let $D$ be a non trivial symmetric $(v,k,\lambda)$ design and
$\alpha$ a non trivial automorphism of $D$. Then $|F(\alpha)|\leq
k+\sqrt{k-\lambda}.$
\end{lem}
\begin{lem}\label{main1}
Let $D$ be a  symmetric $(v,k,\lambda)$ design and $\alpha$ an
automorphism of $D$ of prime order $p$ such that $\lambda<p$. If
$B$ is a block of $D$ such that $|F(\alpha)\cap B|\geq 2$, then
$B^\alpha=B$.
\end{lem}
\begin{proof}
Let $x,y$ be two distinct elements of $F(\alpha)\cap B$. Then
$x,y\in B=B^{\alpha^0}, B^\alpha, \dots, B^{\alpha^{\lambda}}$.
Since every two distinct points are in exactly $\lambda$ blocks,
$B^{\alpha^i}=B^{\alpha^j}$ for some distinct
$i,j\in\{0,1,\dots,\lambda\}$. Thus $B^{\alpha^{i-j}}=B$. Since
$p$ is prime and $1\leq |i-j|\leq \lambda<p$, $\gcd(i-j,p)=1$.
Therefore $B^{\alpha}=B$ as required.
\end{proof}
\begin{lem}\label{main2}
Let $B_1$ and $B_2$ be two distinct fixed blocks of the
automorphism $\alpha$ of prime order $p$ of a symmetric
$(v,k,\lambda)$ design with $\lambda<p$. Then $B_1\cap
B_2\subseteq F(\alpha)$.
\end{lem}
\begin{proof}
Suppose, for a contradiction, that there exists a point $x\in
(B_1\cap B_2)\backslash F(\alpha)$. Thus
$x^{\alpha^i}\not=x^{\alpha^j}$, for any two distinct
$i,j\in\{0,1,\dots,p-1\}$; since otherwise $x^{\alpha^{i-j}}=x$
and so $x^{\alpha}=x$, as $\gcd(i-j,p)=1$. It follows that
$p=|\{x^\beta \;|\; \beta\in\langle \alpha \rangle\}|$. Since
$B_i^\alpha=B_i$ for $i\in\{1,2\}$, we have that  $\{x^\beta
\;|\; \beta\in\langle \alpha \rangle\}\subseteq  B_1\cap B_2$.
Therefore $|B_1\cap B_2|\geq p>\lambda$, a contradiction; since in
symmetric $(v,k,\lambda)$ designs, two distinct blocks intersect
in exactly $\lambda$ points.
\end{proof}
\begin{lem}\label{m4}
Let $\alpha$ be an automorphism  of prime order $p$ of a symmetric
$(v,k,\lambda)$ design with $\lambda<p$. Then
$$|F(\alpha)|+\displaystyle\sum_{B\in F_b(\alpha)}|B\backslash
F(\alpha)|\leq v.$$
\end{lem}
\begin{proof}
It follows from Lemma \ref{main2} that for any two distinct
blocks $B_1$ and $B_2$ in $F_b( \alpha)$, $\big(B_1\backslash
F(\alpha)\big) \cap \big(B_2\backslash
F(\alpha)\big)=\varnothing$. This completes the proof.
\end{proof}
\begin{lem}\label{m3} Let $\alpha$ be an automorphism of a  symmetric $(v,k,\lambda)$
design of prime order $p$ such that $1<\lambda<p$. Then
$B\nsubseteq F(\alpha)$ for all blocks $B$.
\end{lem}
\begin{proof} Suppose, for a contradiction, that there exists a block $B$ such that
$B\subseteq F(\alpha)$.  Since every block $B_1\not=B$ intersects
$B$ in $\lambda\geq 2$ points, it follows from Lemma \ref{main1}
that every block is  fixed under $\alpha$. Thus
$|F_b(\alpha)|=|F(\alpha)|=v$, by Lemma \ref{fix=fix}.  Hence
$\alpha$ is the identity automorphism; a contradiction. This
completes the proof.
\end{proof}
The following lemma is Theorem 2.7 of Aschbacher's paper
\cite{As}.
\begin{lem}{\rm (Theorem 2.7 of \cite{As})}\label{main}
Let $p$  be a prime divisor of the automorphism group of a
symmetric $(v,k,\lambda)$ design such that $1<\lambda<p$ and
$\gcd(p,v)=1$. Then  $p\leq k$.
\end{lem}
\begin{proof}
Suppose that $\alpha$ is an automorphism of the design of order
$p$. Since $\alpha$ is a permutation on the point set,
$F(\alpha)\equiv v \mod{p}$ and since $\gcd(p,v)=1$, we have that
$|F(\alpha)|\geq 1$. Thus, by Lemma \ref{fix=fix}, there exists a
block $B$ such that $B^\alpha=B$. Thus by Lemma \ref{m3}, there
exists an element $x\in B\backslash F(\alpha)$ and so $|\{x^\beta
\;|\; \beta \in \langle \alpha \rangle\}|=p$. Since $B^\alpha=B$,
we have that  $\{x^\beta \;|\; \beta \in \langle \alpha
\rangle\}\subseteq B$ and so $p\leq k$, as required.
\end{proof}
\section{Automorphism group of a possible symmetric (81,16,3) design}
\begin{lem}\label{elem}
Let $G$ be an automorphism group of a possible symmetric
$(81,16,3)$ design which is elementary abelian $2$-group. Then
$|G|\leq 8$.
\end{lem}
\begin{proof}
Let $r$ be the number of orbits of the action of $G$ on the point
set of the design. Then by the Cauchy-Frobenius Lemma (see
\cite[Proposition A.2, p. 246]{Lander}),
$$r=\frac{1}{|G|}\sum_{\alpha\in G} |F(\alpha)|.$$
Since $G$ is an elementary abelian $2$-group, it follows from
Theorem \ref{fix2}, that $|F(\alpha)|=9$ for all non-identity
elements $\alpha$ of $G$. Let $|G|=2^n$. Then,  since
$r=(2^n+8)\cdot 9/2^n$ is an integer, we must have that $2^n$
divides $2^n+8$ and so $n\leq 3$, as required.
\end{proof}
\begin{lem}\label{7,11}
Let $G$ be an automorphism group of a possible symmetric
$(81,16,3)$ design. Then $G$ has no element of order $7$ or $11$.
\end{lem}
\begin{proof}
Suppose, for a contradiction, that $G$ has an automorphism
$\alpha$ of order $p$, where $p\in \{7,11\}$. Since $\alpha$ is a
permutation on a set with $81$ elements,  we have
$|F(\alpha)|\equiv 81  \mod{p}$. Then it follows from Lemma
\ref{fix} that
$$|F(\alpha)|\in\begin{cases} \{4,11,18\}& \; \text{if} \; p=7 \\ \{4,15\} & \; \text{if} \; p=11 \end{cases}. \eqno{(I)}$$
Thus there are at least two distinct blocks which are fixed by
$\alpha$ and so $$|F(\alpha)|\geq 3 \eqno{(*)}$$ by Lemma
\ref{main2}. Now if $B\in F_b(\alpha)$, then $\alpha$ induces a
permutation on the set $B$. Therefore $|F(\alpha) \cap B|\equiv
16 \mod{p}$ and so by $(*)$ and Lemma \ref{m3} we have
$$|F(\alpha) \cap B|=\begin{cases} 9 & \; \text{if} \;
p=7 \\ 5 & \; \text{if}\; p=11
\end{cases}. \eqno{(II)}$$
If $p=11$, then it follows from $(I)$ and $(II)$ that
$|F(\alpha)|=15$ and $|B\backslash F(\alpha)|=11$ for all blocks
$B\in F_b(\alpha)$; and if $p=7$, then   $|B\backslash
F(\alpha)|=7$ for all blocks $B\in F_b(\alpha)$ and
$|F(\alpha)|\in\{11,18\}$. Both cases  contradict Lemma \ref{m4}.
This completes the proof.
\end{proof}
\begin{lem}\label{fix5}
Let $\alpha$ be an automorphism  of a possible symmetric
$(81,16,3)$ design of order $5$. Then $|F(\alpha)|=1$.
\end{lem}
\begin{proof}
Since $\alpha$ is a permutation on the point set, it follows from
Lemma \ref{fix} that  $|F(\alpha)|\in\{1,6,11,16\}$. Suppose, for
a contradiction, that $|F(\alpha)|\not=1$. Let $B=B_1$ be an
arbitrary block in $F_b(\alpha)$. Since
$|F_b(\alpha)|=|F(\alpha)|\geq 2$, there exists a block
$B_2\not=B_1$ in $F_b(\alpha)$. By Lemma \ref{main2}, $B_1\cap
B_2\subseteq F(\alpha)$ and so there exist distinct elements $x$
and $y$ in $F(\alpha)$ which are both in $B_1$ and $B_2$.
Therefore there exists a block  $B_3$ distinct from $B_1$ and
$B_2$  containing both $x$ and $y$. Thus $3=|B_i\cap B_j|\geq
|B_1\cap B_2 \cap B_3|\geq 2$ for any two distinct $i,j\in
\{1,2,3\}$. Now by Lemma \ref{main2}, $B_i^\alpha=B_i$ for all
$i\in\{1,2,3\}$ and so $\alpha$ is a permutation on  $B_i$.
Therefore, it follows from  Lemmas \ref{main2} and \ref{m3}, that
$|F(\alpha)\cap B|\in\{6,11\}$ for all blocks $B\in F_b(\alpha)$.
Thus $|F(\alpha) \cap \big( B_1\cup B_2 \cup B_3\big)|\geq 11$
and so $|F(\alpha)|\in\{11,16\}$. If $|F(\alpha)|=16$, then
$$|F(\alpha)|+\sum_{B\in F_b(\alpha)}|B\backslash F(\alpha)|\geq 16+16\cdot 5=85,$$
which is a contradiction by Lemma \ref{m4}. If $|F(\alpha)|=11$,
then there is no block $B\in F_b(\alpha)$ such that
$|F(\alpha)\cap B|=11$, since otherwise $|(B'\cup B) \cap
F(\alpha)|\geq 11+6-3=14$ for any block $B'\in F(\alpha)$
distinct from $B$. Hence, in this case, $$|F(\alpha)|+\sum_{B\in
F_b(\alpha)}|B\backslash F(\alpha)|\geq 11+11\cdot 10=121,$$
which contradicts Lemma \ref{m4}. This completes the proof.
\end{proof}
\begin{lem}\label{5}
Let $G$ be an automorphism group of a possible symmetric
$(81,16,3)$ design which is a $5$-group. Then $|G|\leq 5$.
\end{lem}
\begin{proof}
It is enough to show that $G$ has no subgroup $H$ of order $5^2$.
If $\alpha\in G$ is of order $25$, then by Lemma \ref{fix5},
$|F(\alpha)|=1$, since $\varnothing\not=F(\alpha)\subseteq
F(\alpha^5)$.
   Then,
by Lemma \ref{fix5}, the number of orbits of the action of $H$ on
$G$ is equal to
$$r=\frac{1}{5^2}\sum_{h\in H}|F(h)|=\frac{81+24\cdot 1}{5^2}=\frac{21}{5}.$$
This is a contradiction, since $r$ should be an integer.
\end{proof}
\begin{lem}\label{fix13}
Let $\alpha$ be an automorphism  of a possible symmetric
$(81,16,3)$ design of order $13$. Then $|F(\alpha)|=3$.
\end{lem}
\begin{proof}
Since $\alpha$ is a permutation on the point set, it follows from
Lemma \ref{fix} that  $|F(\alpha)|\in\{3,16\}$. Suppose, for a
contradiction, that $|F(\alpha)|=|F_b(\alpha)|=16$. Then, by Lemma
\ref{m3},  $|F(\alpha)\cap B|=3$ for all $B\in F_b(\alpha)$. Thus
$$|F(\alpha)|+\sum_{B\in F_b(\alpha)}|B\backslash F(\alpha)|\geq
16+ 16 \cdot 13=224,$$ contradicting Lemma \ref{m4}. This
completes the proof.
\end{proof}
\begin{lem}\label{13}
Let $G$ be an automorphism group of a possible symmetric
$(81,16,3)$ design which is a $13$-group. Then $|G|\leq 13$.
\end{lem}
\begin{proof}
It is enough to show that $G$ has no subgroup $H$ of order $13^2$.
Since $13^2>81$, $G$ has no element of order $13^2$. Thus $H$ is
an elementary abelian $13$-group. Then, by Lemma \ref{fix13}, the
number of orbits of the action of $H$ on $G$ is equal to
$$r=\frac{1}{13^2}\sum_{h\in H}|F(h)|=\frac{81+12\cdot 3}{13^2}=\frac{9}{13}.$$
This is a contradiction, since $r$ should be an integer.
\end{proof}
\begin{lem}\label{el}
Let $G$ be an automorphism group of a possible symmetric
$(81,16,3)$ design. Then $G$ has no element with the following
orders: $10$, $26$,  $65$.
\end{lem}
\begin{proof}
\begin{enumerate}
\item  Suppose that $G$ has an element of order $10$. Then $G$
contains two automorphisms $\alpha$ and $\beta$ of orders $5$ and
$2$ respectively such that $\alpha \beta=\beta \alpha$. Since
$\alpha$ and $\beta$ commutes,
$\alpha\big(F(\beta)\big)=F(\beta)$. By Theorem \ref{fix2}  we
have that  $|F(\beta)|=9$.  Now by considering the cycle
decomposition of $\alpha$ on $F(\beta)$, it follows that
$|F(\alpha)\cap F(\beta)|\in\{4,9\}$ which contradicts Lemma
\ref{fix5}.

\item Suppose that $G$ has an element of order 26. Then  $G$
contains two automorphisms $\alpha$ and $\beta$ of orders $13$
and $2$ respectively such that $\alpha \beta=\beta \alpha$. Since
$\alpha$ and $\beta$ commutes,
$\alpha\big(F(\beta)\big)=F(\beta)$ and  by Theorem \ref{fix2},
$|F(\beta)|=9$,   the cycle decomposition of $\alpha$ on
$F(\beta)$ shows that $F(\beta)\subseteq F(\alpha)$ which
contradicts Lemma \ref{fix13}.

\item Suppose that $G$ has an element of order 65. Then  $G$
contains two automorphisms $\alpha$ and $\beta$ of orders $13$
and $5$ respectively such that $\alpha \beta=\beta \alpha$. Since
$\alpha$ and $\beta$ commutes,
$\beta\big(F(\alpha)\big)=F(\alpha)$. But by  Lemma \ref{fix13} we
have that $|F(\alpha)|=3$ so the cycle decomposition of $\beta$
on $F(\alpha)$ implies that $F(\alpha)\subseteq F(\beta)$ which
contradicts Lemma \ref{fix5}.
\end{enumerate}
\end{proof}
\noindent{\bf Proof of Theorem \ref{mainmain}.} It follows from
Lemmas \ref{elem}, \ref{7,11}, \ref{5}, \ref{13} and \ref{el}

\end{document}